\documentclass{article}
\usepackage{amsmath}

\input{tcilatex}
\begin{document}

\begin{center}
GAUSSIAN ESTIMATES: A BRIEF HISTORY

\bigskip
\end{center}

\begin{flushleft}
D.G.Aronson

School of Mathematics

University of Minnesota

Minneapolis MN 55455

e-mail: arons001@umn.edu

\bigskip

ABSTRACT: Two-sided Gaussian estimates for the fundamental solution of a
second order linear parabolic differential equation are upper and lower
bounds in terms of the fundamental solution of the classical heat conduction
equation. In his seminal 1958 paper Nash stated, without proof, two-sided
non-Gaussian bounds for the fundamental solution of a uniformly parabolic
divergence structure equation assuming only boundedness of the coefficients.
In his 1967-1968 papers Aronson derived truly Gaussian estimates for the
fundamental solutions of a large class of linear parabolic equations
(including the divergence structure equations) under minimal non-regularity
assumptions on the coefficients. Subsequently in 1986 Fabes \& Stroock
derived Gaussian estimates for the divergence structure equation directly
from the ideas of Nash and went on to prove Nash's continuity theorem and
the Harnack inequality as a consequence of their estimate. In this note I
describe these results together with various extensions.

\bigskip
\end{flushleft}

Let $\Gamma (x,t;\xi ,\tau )$ denote the\textit{\ fundamental solution} of
the divergence structure second order linear parabolic equation%
\begin{equation}
\partial _{t}u-\tsum_{j=1}^{N}\partial _{x_{J}}\left\{
\tsum_{i=1}^{N}a_{ij}(x,t)\partial _{x_{i}}u+a_{j}(x,t)u\right\}
-\tsum_{j=1}^{N}b_{j}(x,t)\partial _{x_{j}}u-c(x,t)u=0,  \tag{1}
\end{equation}%
where $(x,t)\in S=\mathbf{R}^{N}\times (0,T)$ for some $T>0.$ If the
coefficients of (1) are smooth then the fundamental solution exists in the
classical sense, otherwise it must be interpreted in the weak sense (cf.
[2]). The \textit{two-sided Gaussian estimate} for $\Gamma $ is the
following:

\textbf{Theorem: }\textit{There exist positive constants }$\alpha
_{1},\alpha _{2}$\textit{\ and }$C$\textit{\ depending only on }$T$\textit{\
and the structure of equation (1) such that}%
\begin{equation}
C^{-1}g_{1}(x-\xi ,t-\tau )\leq \Gamma (x,t:\xi ,\tau )\leq Cg_{2}(x-\xi
,t-\tau )  \tag{2}
\end{equation}%
\textit{for all }$(x,t),(\xi ,\tau )\in S$\textit{\ with }$t>\tau $\textit{,
where }$g_{i}(x,t)$\textit{\ denotes the (Gaussian) fundamental solution of
the heat equation }$\alpha _{i}\Delta u=\partial _{t}u$\textit{\ for }$%
i=1,2. $

\begin{flushleft}
Note that this estimate does not require any smoothness assumption on the
coefficients of equation (1).
\end{flushleft}

The Gaussian estimate (2) was first proved in 1967 for the special case of
the equation%
\begin{equation}
\partial _{t}u-\tsum_{j=1}^{N}\partial _{x_{j}}\left\{
\tsum_{i=1}^{N}a_{ij}(x,t,)\partial _{x_{i}}u\right\} =0  \tag{3}
\end{equation}
in [1] and subsequently extended to the general case of equation (1) in [2].
In his seminal 1958 paper [10] John Nash considers equation (3) and proves H%
\"{o}lder continuity of solutions assuming only uniform parabolicity and
bounded measurable coefficients (as is also done in [1]). In the appendix to
his paper he states two-sided estimates for the fundamental solution, but
his estimates are not Gaussian and complete proofs are not provided. He also
points out that it is possible to derive the Harnack inequality from his
bounds for the fundamental solution.

There are two distinct methods of deriving Gaussian estimates for the
fundamental solution of (3). In the original 1967 derivation [1] the
estimates are obtained as a consequence of an energy estimate, along with
the Harnack inequality and H\"{o}lder continuity of solutions proved in [9]
and [3]. On the other hand in their 1986 paper, Fabes\ \& Stroock made a
detailed study of Nash's work and discovered that his method can be made to
yield the two-sided Gaussian estimate directly. They then showed that the
estimate can be used to derive the H\"{o}lder continuity of solutions and
the Harnack inequality. Neither Nash or Fabes \& Stroock consider the full
equation (1). Thus the Gaussian estimates for the full equation (1) are
still only known as a consequence of the various properties of weak
solutions proved in [9] and [3].

In [1] it is assumed that the coefficients $a_{ij}(x,t)$ in equation (3) are
bounded and measurable in $S$, and that there exists a constant $\nu >0$
such that (with summation over repeated indices)%
\begin{equation*}
a_{ij}(x,t)\zeta _{i}\zeta _{j}\geq \nu \left\vert \zeta \right\vert ^{2}
\end{equation*}%
almost everywhere in $S$ for all $\zeta \in \mathbf{R}^{N}.$ The existence
of the weak fundamental solution $\Gamma $ is proved under conditions which
include these in [2]. The constants $\alpha _{1},\alpha _{2}$ and $C$ in the
estimate (2) depend only on $\nu ,N,T$ and the bounds for the coefficients
of (3). For the general equation (1) the Gaussian bounds are proved in [2]
under the following assumptions on the structure of equation (3) which will
be referred to collectively as (H). There exist constants $0<\nu ,M<\infty $
and $0\leq M_{0}<\infty $ such that the coefficients of equation (1) satisfy

(H.1) For all $\zeta \in \mathbf{R}^{N}$ and for almost all $(x,t)\in S$%
\begin{equation*}
a_{ij}(x,t)\zeta _{i}\zeta _{j}\geq \nu \left\vert \zeta \right\vert ^{2}%
\text{ and }\left\vert a_{ij}\right\vert \leq M.
\end{equation*}

(H.2) Each of the coefficient $a_{j}$ and $b_{j}$ belong to some Bochner
space $L^{pq}(S)$, where%
\begin{equation*}
2<p,q\leq \infty \text{ and }\frac{N}{2p}+\frac{1}{q}<\frac{1}{2},
\end{equation*}%
and $\left\vert a_{j}(x,t\right\vert ,\left\vert b_{j}(x,t\right\vert \leq
M_{0}$ almost everywhere in $S.$

(H.3) $c\in L^{pq}(S),$ where%
\begin{equation*}
1<p,q\leq \infty \text{ and }\frac{N}{2p}+\frac{1}{q}<1
\end{equation*}%
and $c(x,t)\leq M_{0}$ almost everywhere in $S.$ \ \ \ \ \ \ \ \ \ \ \ \ \ \
\ \ \ \ \ \ \ \ \ \ \ \ \ \ \ \ \ \ \ \ \ \ \ \ \ \ \ \ \ \ \ \ \ \ \ \ \ \
\ \ \ \ \ \ \ \ \ \ \ \ \ \ \ \ \ \ \ \ \ \ \ \ \ \ \ \ \ \ \ \ \ \ \ \ \ \
\ \ \ \ \ \ \ \ \ \ \ \ \ \ \ \ \ \ \ \ \ \ \ \ \ \ \ \ \ \ \ \ \ \ \ \ \ \
\ \ \ \ \ \ \ \ \ \ \ \ \ \ \ \ \ \ \ \ \ \ \ \ \ \ \ \ \ \ \ \ \ \ \ \ \ \
\ \ \ \ \ \ \ \ \ \ \ \ \ \ \ \ \ \ \ \ \ \ \ \ \ \ \ \ \ \ \ \ \ \ \ \ \ \
\ \ \ \ \ \ \ \ \ \ \ \ \ \ \ \ \ \ \ \ \ \ \ \ \ \ \ \ \ \ \ \ \ \ \ \ \ \
\ \ \ \ \ \ \ \ The existence of a weak fundamental solution under the
hypothesis (H) is proved in [2]. The constants $\alpha _{1},\alpha _{2}$ and 
$C$ in the estimate (2) depend only on $N,$ $T$ and the structure of
equation (1), i.e., the hypothesis (H).

The proof of the lower bound in (2) in both \ [1] and [2] rely heavily on
consequences of the regularity of solutions and the pointwise Harnack
inequality proved in [3] and [9]. On the other hand the proof of the upper
bound does not involve the Harnack inequality or its consequences. Instead
it uses the estimate (first proved by Nash in the special case of equation
(3))%
\begin{equation}
\Gamma (x,t;\xi ,\tau )\leq C(t-\tau )^{-N/2}  \tag{4}
\end{equation}%
in $S\times S$ for $t>\tau $, where the constant $\mathcal{C}$ depends only
on $N,T$ and (H); a technical lemma which gives an estimate for the
magnitude inside a ball of a solution which is initially supported in the
exterior of that ball; and the semigroup (reproducing) property of the
fundamental solution.

\ \ \ In his book [4], Davis considers equation (3) in the case the
coefficients $a_{ij}$ are independent of $t$ and depend only on $x.$ He
introduces a technique which enables him to sharpen the upper bound in (2)
by replacing the Euclidian distance with the Riemannian distance associated
with the coefficients. Fabes \& Stroock [6] refine Nash's sketch of the
argument he indicated to establish his non-Gaussian upper bound and use it
together with Davis' method to prove the Gaussian upper bound (2) for the
fundamental solution of equation (3). They point out \textquotedblleft 
\textit{that the upper bound itself is an important tool for our
understanding and simplification of those ideas of Nash needed to obtain the
lower bound'}" in (2). \ Although their procedure is basically due to Nash,
the Gaussian upper bound allows them to simplify his argument and refine his
lower bound. The derivation of the lower bound in [6] depends on two
estimates,.both essentially due to Nash: the inequality (4) and the
existence of a constant $B<\infty $ depending only on $\nu $ such that for
all%
\begin{equation*}
\left\vert x\right\vert \leq 1\tint e^{-\pi \left\vert y\right\vert
^{2}}\log \Gamma (x,1;y,0)dy\geq -B.
\end{equation*}%
Using these estimates together with the semigroup property of the
fundamental solution yield the Gaussian lower bound. Using the Gaussian
bound (2), Fabes \& Stroock go on to derive Nash's H\"{o}lder continuity
result and the Harnack inequality. For the latter result they use the
methods of Krylov \& Safonov [7]. It should be emphasized that Fabes \&
Stroock derive their bounds directly from equation (3) without any reference
to regularity properties of the solution and, indeed derive the regularity
results as a consequence of their bounds.\ Fabes [5] extended Nash's
`moment' bound to estimate all of the moments of the fundamental solution of
(3) and applies Davis' techniques to prove Davis' Riemannian upper bound. He
also shows that these methods can yield upper bounds for heat kernels on a
class of complete Riemann manifolds.

If the coefficients $a_{ij}$ in equation (3) are independent of $t$ and $%
N\geq 3,$ then 
\begin{equation*}
\int_{0}^{\infty }\Gamma (x,t;\xi ,0)=G(x,\xi ),
\end{equation*}%
where $G(x,\xi )$ is the fundamental solution of the elliptic equation%
\begin{equation*}
\tsum_{j=1}^{N}\partial _{x_{j}}\tsum_{i=1}^{N}a_{ij}(x)\partial _{x_{i}}u=0.
\end{equation*}%
As noted in [1], in this case the constants in estimate (2) can be chosen
independent of $T$ and we can integrate these estimates over $R^{+}$ to
obtain the estimate%
\begin{equation*}
K^{-1}\left\vert x-\xi \right\vert ^{2-N}\leq G(x,\xi )\leq K\left\vert
x-\xi \right\vert ^{2-N}.
\end{equation*}%
This estimate is not new having been derived directly from potential
theoretic considerations by Littman, Stampacchia \& Weinberger [8] and H.
Royden [13].

Porper \& Eidel'man [12] consider a slight generalization of equations (1)
and (3) involving a coefficient $p(x)$ multiplying $\partial _{t}u$. Using
essentially the arguments employed in [1] and [2], they give a detailed
account of the derivation of H\"{o}lder continuity, the Harnack inequality
and the two-sided Gaussian estimate for the analogue of equation (3) and a
very brief description of these results of the analogue of equation (1)
under conditions similar to (H).

In [11] Norris and Stroock consider the operator%
\begin{equation*}
L\equiv \nabla \cdot (A(x,t)\nabla +AE(x,t)\cdot \nabla -\nabla \cdot (A\hat{%
E}(x.t))+C(x,t),
\end{equation*}%
where the coefficients of $L$ are measurable functions on $\mathbf{R}%
^{N}\times \mathbf{R}$. Here $A$ is an $N\times N$ positive-definite
symmetric matrix, $E$ and $\hat{E}$ are in $\mathbf{R}^{N}$ and \ $C$ is in $%
\mathbf{R.}$ Their main result is a very precise two-sided estimate for the
fundamental solution of $Lu=\partial _{t}u$ based on energy functions
associated with the coefficients of $L.$ However they are forced to assume
the uniform continuity of $A$ and $E-\hat{E}.$

\bigskip

\begin{flushleft}
\textbf{REFERENCES}
\end{flushleft}

1. \ D.G.Aronson \textit{Bounds for the Fundamental Solution of a Parabolic
Equation,}\ Bull. AMS \textbf{73}(1967), 890-896.

2. \ D.G.Aronson \textit{Non-Negative Solutions of Linear Parabolic Equations%
}, Ann. Scuola Norm. Sup. Pisa \textbf{22}(1968), 607-694: Non-Negative%
\textit{\ Solutions of Linear Parabolic Equations: An Addendum,} Ann. Scuola
Norm. Sup. Pisa \textbf{25}(1971), 221-228.

3. \ D.G.Aronson \& J.Serrin, \textit{Local Behavior of Solutions of
Quasilinear Parabolic Equations}, Arch. Rat. Mech. Anal. \textbf{25}(1967),
81-122.

4. \ E.B.Davis\textbf{\ Heat Kernels and Spectral Theory}, Cambridge Tracts
in Mathematics \textbf{92, }Cambridge University Press 1989.

5. \ E.B.Fabes \textit{Gaussian Upper Bounds on the Fundamental Solution of
Parabolic Equations: the Method of Nash}, Lecture Notes in Mathematics 
\textbf{1563}(1993), 1-20, Springer-Verlag.

6. \ E.B.Fabes \& D.W.Stroock\textit{\ A New Proof of Moser's Parabolic
Harnack Inequality Using Old Ideas of Nash}, Arch. Rat. Mech. Anal. \textbf{%
96}(1986), 327-338.

7. \ N.V.Krylov \& M.V.Safonov \textit{A certain property of solutions of
parabolic equations with measurable coefficients}, Math. USSR Izv. \textbf{16%
}(1981), 151-164.

8. \ W.Littman, G.Stampacchia \& H.F.Weinberger \textit{Regular Points for
Elliptic Equations with Discontinuous Coefficients}, Ann. Scuola Norm. Sup.
Pisa \textbf{17}(1968), 43-77.

9. \ J.Moser\textit{\ A Harnack Inequality for Parabolic Differential
Equations} Comm. Pure Appl. Math. \textbf{17}(1964), 101-134; \textit{A
Correction to \textquotedblleft A Harnack Inequality for Parabolic
Differential Equations"}, Comm. Pure Appl. Math. \textbf{20}(1967), 231-236.

10. J.Nash\textit{\ Continuity of Solutions of Parabolic and Elliptic
Equations}, Amer. J. Math. \textbf{80}(1958), 931-954.

11. James R. Norris \& Daniel W. Stroock \textit{Estimates on the
Fundamental Solution to Heat Flows with Uniformly Elliptic Coefficients},
Proc. London Math. Soc. \textbf{62}(1991), 373-402.

12. \ F.O.Proper \& S.D.Eidel'man \textit{Two-Sided Estimates of the
Fundamental Solution of Second Order Parabolic Equations and Some
Applications,} Russian Math. Surveys \textbf{39}(1984), 119-178.

13. H.Royden \textit{The Growth of a Fundamental Solution of an Elliptic
Divergence Structure Equation}, \textbf{Studies in Mathematical Analysis and
Related Topics: Essays in Honor of George Polya, }Stanford 1962, 333-340.

\end{document}